\documentclass[11pt]{article}

\usepackage{amsmath,amssymb}

\newtheorem{as}{Theorem}[section]
\newtheorem{theorem}[as]{Theorem}
\newtheorem{prop}[as]{Proposition}
\newtheorem{lemma}[as]{Lemma}

\newtheorem{conj}[as]{Conjecture}

\newtheorem{definition}[as]{Definition}

\newcommand{\qed}{\hspace*{\fill} \rule{7pt}{7pt}}
\newcommand{\Proof}{\noindent{\bf Proof.}\ \ }

\begin{document}

\title{Acyclic subgraphs with high chromatic number \thanks{This research was supported by the Israel Science Foundation (grant No. 1082/16).}}

\author{
Safwat Nassar
\thanks{Department of Mathematics, University of Haifa, Haifa
31905, Israel.}
\and
Raphael Yuster
\thanks{Department of Mathematics, University of Haifa, Haifa
31905, Israel. Email: raphy@math.haifa.ac.il}
}

\date{}

\maketitle

\setcounter{page}{1}

\begin{abstract}

For an oriented graph $G$, let $f(G)$ denote the maximum chromatic number of an acyclic subgraph of $G$.
Let $f(n)$ be the smallest integer such that every oriented graph $G$ with chromatic number larger than $f(n)$ has $f(G) > n$.
Let $g(n)$ be the smallest integer such that every tournament $G$ with more than $g(n)$ vertices has $f(G) > n$.
It is straightforward that $\Omega(n) \le g(n) \le f(n) \le n^2$.
This paper provides the first nontrivial lower and upper bounds for $g(n)$.
In particular, it is proved that $\frac{1}{4}n^{8/7} \le g(n) \le n^2-(2-\frac{1}{\sqrt{2}})n+2$.
It is also shown that $f(2)=3$,  i.e. every orientation of a $4$-chromatic graph has a $3$-chromatic
acyclic subgraph.
Finally, it is shown that a random tournament $G$ with $n$ vertices has $f(G) = \Theta(\frac{n}{\log n})$ whp.
\end{abstract}

\section{Introduction}

All graphs in this paper are finite and simple. An {\em orientation} of an undirected graph is obtained by assigning a direction to each
edge. An important class of oriented graphs are {\em tournaments} which are
orientations of a complete graph. We denote by $T_n$ the unique acyclic (thereby transitive) tournament with $n$ vertices.
An acyclic subgraph of an oriented graph is a subgraph having no directed cycles.
In this paper, the chromatic number of an oriented graph is the chromatic number of its underlying undirected graph.

It is a folklore argument that every oriented graph has acyclic subgraphs containing at least half of the edges. Indeed,
every linear ordering of the vertices partitions the edge set to two acyclic subgraphs, one consisting of the edges pointing from lower vertices to
higher vertices and the other consisting of the edges pointing from higher vertices to lower vertices. At least one of these two subgraphs
contains at least half of the edges. Thus, we are guaranteed to find dense (with respect to the density of the original oriented graph) acyclic 
subgraphs. However, apart from density, we do not have much control on other complexity parameters of these acyclic subgraphs.
Perhaps the most natural is the following question raised by Addario-Berry et al. \cite{AHSRT-2013} and also in the survey of Havet
\cite{havet-preprint}. Suppose we know that the chromatic number of our oriented graph is large, can we guarantee that an acyclic subgraph of
it also has high chromatic number?

\begin{definition}\label{def:1}
For an oriented graph $G$, let $f(G)$ denote the maximum chromatic number of an acyclic subgraph of $G$.
Let $f(n)$ be the smallest integer such that every oriented graph $G$ with chromatic number larger than $f(n)$ has $f(G) > n$.
\end{definition}
It was observed in \cite{AHSRT-2013} that $f(n) \le n^2$ by the following standard ``product-coloring'' argument.
Suppose $G$ is a graph, and take any partition of its edge set into $k$ parts, inducing subgraphs $G_1, \ldots, G_k$.
Then clearly, $\prod_{i=1}^k \chi(G_i) \ge \chi(G)$ as we may properly color each vertex of $G$ by the $k$-dimensional vector whose $i$'th entry is
the color that vertex received in a coloring of $G_i$ with $\chi(G_i)$ colors. In particular, for the case where $G$ is an oriented graph with
$\chi(G) \ge n^2+1$ and $G_1$ and $G_2$ are acyclic subgraphs forming a partition of the edge set according to some linear ordering of the vertices
as described in the previous paragraph,
then $\chi(G_1) \cdot \chi(G_2) \ge n^2+1$ so at least one of them has chromatic number at least $n+1$.
This shows that $f(n) \le n^2$. To date, there is no known improvement over this simple upper bound.

What if we restrict the question to tournaments? A tournament with $g(n)$ vertices has chromatic number $g(n)$ as it is an orientation of $K_{g(n)}$.
Thus, we have the following definition.
\begin{definition}\label{def:2}
Let $g(n)$ be the smallest integer such that every tournament $G$ with more than $g(n)$ vertices has $f(G) > n$.
\end{definition}
Clearly, we have $g(n) \le f(n)$ and hence the aforementioned simple upper bound $g(n) \le n^2$ holds here as well.
Our first result is a modest, yet nontrivial improvement to the upper bound.
\begin{theorem}\label{t:1}
$g(n) \le n^2-(2-\frac{1}{\sqrt{2}})n+2$.
\end{theorem}
We suspect that this upper bound can be improved and raise the following conjecture.
\begin{conj}\label{conj:1}
$g(n) = o(n^2)$.
\end{conj}
As for $f(n)$, while we cannot improve upon the upper bound $f(n) \le n^2$ in general,
we do settle the first nontrivial case.
\begin{theorem}
\label{t:3}
Suppose $G$ is an orientation of a $4$-chromatic graph. Then $G$ has an acyclic subgraph with chromatic
number at least $3$. In particular, $f(2)=3$.
\end{theorem}
Notice that $g(2) = 3$ is trivial since every tournament with more than $3$ vertices has an acyclic triangle
and since each acyclic subgraph of the directed triangle has chromatic number at most $2$.
Whether $f(n)=g(n)$ for larger $n$ remains open.

Let $h(n)$ be the least integer such that every oriented graph with chromatic number $h(n)$ contains every oriented tree
with $n$ vertices. A longstanding conjecture of Burr \cite{burr-1980} asserts that $h(n)=2n-2$.
Presently the best known upper bound is $h(n) \le n^2/2-n/2+1$ given in \cite{AHSRT-2013}.
One motivation for studying $f(n)$, given in \cite{AHSRT-2013}, stems from the fact proved there that every acyclic oriented graph with chromatic
number $n$ contains every oriented tree on $n$ vertices. Hence, any upper bound for $f(n)$ can be used as an upper bound for $h(n)$ and any
lower bound for $f(n)$ which is significantly larger than linear shows that obtaining improvements for $h(n)$ are limited with this approach.
Our next result provides such a lower bound for $g(n)$ and hence a lower bound for $f(n)$.
\begin{theorem}\label{t:2}
There are tournaments $G$ with more than $n^{8/7}/4$ vertices such that $f(G) \le n$.
Consequently, $g(n) \ge n^{8/7}/4$.
\end{theorem}
In spite of the fact proved in Theorem \ref{t:2} that there are tournaments that have the property that
every acyclic subgraph has chromatic number polynomially smaller than the number of vertices of the tournament,
it turns out that almost all tournaments do have acyclic subgraphs with only a logarithmic fraction loss in the chromatic number.
Recall that the random tournament ${\cal G}(n)$ is the probability space of $n$-vertex tournaments where the orientation
of each edge is decided independently by a fair coin flip. The value of $f(G)$ for a randomly sampled element of ${\cal G}(n)$
is therefore a random variable, and saying that with high probability $f(G) = \Theta(n / \log n)$ means that there are absolute constants
$c$ and $C$ such that asymptotically almost surely (i.e. with probability tending to $1$ as $n$ tends to infinity)
$cn/\log n \le f(G) \le Cn/\log n$.
\begin{prop}\label{prop:1}
For a random tournament $G$ on $n$ vertices, with high probability $f(G) = \Theta(\frac{n}{\log n})$.
\end{prop}
In fact, the proof of Proposition \ref{prop:1} shows that $f(G) = (1+o_n(1))\frac{n}{2\log n}$ asymptotically almost surely
(all logarithms hereafter are in base $2$).
Random tournaments have the property that their largest transitive subtournaments have logarithmic size.
On the other hand, very ``non-random'' tournaments have transitive subtournaments that are very large,
so one might suspect that random tournaments yield constructions that are not far from $g(n)$.
The proof of Theorem \ref{t:2} show that this is not the case.

Hereafter, the proof of Theorem \ref{t:1} upper bounding $g(n)$ appears in Section 2.
The proof of Theorem \ref{t:2} lower bounding $g(n)$ is given in Section 3.
Section 4 consists of the proof of Proposition \ref{prop:1} concerning $f(G)$ of random tournaments.
Section 5 contains the proof of Theorem \ref{t:3}.

\section{Upper bound for $g(n)$}

Let $G$ be a tournament. For a linear ordering $\pi$ of $V(G)$, consider the ordering of the vertices of $G$ according to
$\pi$. Let $R_\pi(G)$ be the spanning subgraph of $G$ where
an edge $(i,j)$ of $G$ is an edge of $R_\pi(G)$ if and only if $i$ appears before $j$ in the ordering (thus, all edges of $R_\pi(G)$ go from left to right).
Similarly, let $L_\pi(G)$ be the spanning subgraph of $G$ where
an edge $(i,j)$ of $G$ is an edge of $L_\pi(G)$ if $i$ appears after $j$ in the ordering (thus, all edges of $L_\pi(G)$ go from right to left).
Both $R_\pi(G)$ and $L_\pi(G)$ are therefore acyclic and together form a $2$-partition the edge set of $G$.
Also notice that $R_\pi(G)=L_\sigma(G)$ where $\sigma$ is the reverse of $\pi$.

\begin{lemma}\label{l:1}
Suppose $G$ is a tournament with $t$ vertices. Then there is a linear ordering $\pi$ of $V(G)$ such that $R_\pi(G)$ has no transitive sub-tournament
of order $s = \lfloor \sqrt{2t} \rfloor +1$.
\end{lemma}
\Proof
We construct the ordering $\pi$ according to the following greedy procedure. Let $W_1$ be the vertex set of some $T_s$-subgraph of $G$
(if no such subgraph exists then the lemma trivially holds). Suppose the vertices of $W_1$ are $v_{1,1},\ldots,v_{1,s}$
where $(v_{1,j},v_{1,i}) \in E(G)$ if $1 \le i < j \le s$. We start the linear ordering $\pi$ with $v_{1,1},\ldots,v_{1,s}$.
Observe that $W_1$ induces
an independent set in $R_\pi(G)$ (regardless of how we will later extend $\pi$ to the remaining vertices).
Also observe that any $T_s$-subgraph of $G$ that intersects $W_1$ in more than one vertex will not be in $R_\pi(G)$, so it is
{\em destroyed}.

Now, if all $T_s$-subgraphs are destroyed, we are done (we can extend $\pi$ to the remaining $t-s$ vertices
arbitrarily). So, assume that not all $T_s$-subgraphs and destroyed, and let $X$ be the vertex set of some non-destroyed $T_s$-subgraph of $G$.
Since $|W_1 \cap X| \le 1$, we may set $W_2=\{v_{2,1},\ldots,v_{2,s-1}\} \subseteq X \setminus W_1$ and we assume that
$(v_{2,j},v_{2,i}) \in E(G)$ if $1 \le i < j  \le s-1$. We extend $\pi$ with $v_{2,1},\ldots,v_{2,s-1}$. Observe that any $T_s$-subgraph
that intersects $W_2$ in more than one vertex will not be in $R_\pi(G)$, so it is now destroyed. Moreover, any non-destroyed
$T_s$-subgraph has at most one vertex in $W_1$ and and most one vertex in $W_2$.

We continue this process, defining $W_3,W_4,$ etc.
where when we reach stage $k$, we either have all $T_s$-subgraphs destroyed, or we have a set $W_k=\{v_{k,1},\ldots,v_{k,s-k+1}\}$
where we assume that $(v_{k,j},v_{k,i}) \in E(G)$ if $1 \le i < j  \le s-k+1$ and we extend $\pi$ with 
$v_{k,1},\ldots,v_{k,s-k+1}$. Furthermore, we have the property that any $T_s$-subgraph that is not destroyed by now has at most
one vertex in each of $W_1,\ldots,W_k$.

We observe that after $k$ stages, $\pi$ already contains
$\min\{t,\sum_{i=1}^k s-i+1\}  = \min\{t,ks-\binom{k}{2}\}$ vertices, and since $s=\lfloor \sqrt{2t} \rfloor +1$ and the overall number of vertices
is $t$ we must have that after at most $s-1$ stages we have completed the linear order.
Now, since there are only $s-1$ stages, each $T_s$-subgraph of $G$ must be destroyed since it cannot contain two vertices from some
$W_i$. 
\qed

\vspace*{10pt}
We require also the following lemma regarding the minimum number of cyclic triangles that contain an edge of a tournament.
\begin{lemma}\label{l:2}
Suppose $G$ is a tournament with $t$ vertices. Then there is an edge of $G$ that is contained 
in at most $(t+1)/4$ cyclic triangles.
\end{lemma}
\Proof
It is well-known \cite{goodman-1959} that the maximum number of cyclic triangles in a tournament is at most $(t+1)t(t-1)/24$.
The overall number of edges in all of these cyclic triangles is thus at most $(t+1)t(t-1)/8$ and as there are $t(t-1)/2$
edges, some edge appears on at most $(t+1)/4$ cyclic triangles.
\qed

\vspace*{10pt}
\noindent
{\bf Proof of Theorem \ref{t:1}}.
Let $G$ be a tournament on $m= n^2-\lceil (2-\frac{1}{\sqrt{2}})n \rceil +3$
vertices.
We will prove that it has an acyclic subgraph which is not $n$-chromatic.
This will show that $g(n) \le n^2-(2-\frac{1}{\sqrt{2}})n +2$.
Let $(u,v)$ be an edge that is contained in $q \le (m+1)/4$ cyclic triangles. By Lemma \ref{l:2}, such an edge exists.
Let $N^+(u)$ (resp. $N^{-}(u)$) be the set of out-neighbors
(resp. in-neighbors) of $u$ and hence $v \in N^+(u)$. Let $Q$ be the set of vertices such that for each $w \in Q$, the $3$-set $\{u,v,w\}$ forms
a cyclic triangle. Thus, $Q \subseteq N^{-}(u)$ and $|Q|=q \le (m+1)/4$.

Let $G[Q]$ be the induced sub-tournament of $G$ on the set $Q$. Let $\sigma$ be an ordering of $Q$ such that
$R_\sigma(G[Q])$ has no transitive sub-tournament of order $s = \lfloor \sqrt{2q} \rfloor +1$.
By Lemma \ref{l:1}, $\sigma$ exists.
Let $X = N^+(u) \cap N^+(v)$ and let $Y=N^+(u) \cap N^{-}(v)$. We notice that $X \cup Y \cup \{v\} = N^+(u)$.

We construct a linear order $\pi$ of $V(G)$ which consists of the following $6$ consecutive parts.
The prefix of $\pi$ is $X$ (the internal ordering of $X$ is arbitrary).
Following that, we place $v$.
Following that, we place $Y$ (the internal ordering of $Y$ is arbitrary).
Next, we place $u$. 
Next, we place the vertices of $Q$ according to the ordering $\sigma$.
Finally, we place the remaining vertices $N^-(u) \setminus Q$ as the suffix of $\pi$ where the internal ordering of
the vertices of $N^-(u) \setminus Q$ is arbitrary.

We will prove that at least one of $R_\pi(G)$ or $L_\pi(G)$ is not $n$-chromatic.
If one of $R_\pi(G)$ or $L_\pi(G)$ contains a $T_{n+1}$, then we  are done since such a $T_{n+1}$, being a clique, is not $n$-chromatic.
So we may assume that the maximum size of a transitive subtournament in $R_\pi(G)$ and in $L_\pi(G)$ is at most $n$.
Assume, in contradiction, that both $R_\pi(G)$ and $L_\pi(G)$ are $n$-chromatic.
Hence let $C_R=\{R_1,\ldots,R_n\}$ be an $n$-coloring of $R_\pi(G)$ where the $R_i$ are the color classes
and let $C_L=\{L_1,\ldots,L_n\}$  be an $n$-coloring of $L_\pi(G)$ where the $L_i$ are the color classes.
So, each $R_i$ is an independent set in $R_\pi(G)$ and a clique (i.e. a transitive tournament) in $L_\pi(G)$.
Similarly, each $L_i$ is an independent set in $L_\pi(G)$ and a clique in $R_\pi(G)$.
Since there are no cliques of order $n+1$ in each of them, we have that $0 \le |R_i| \le n$ and $0 \le |L_i| \le n$
for $i=1,\ldots,n$.

Now, consider the vertex $u$ in the coloring of $L_\pi(G)$. By the construction of $\pi$, $u$ is connected
in $L_\pi(G)$ to all other vertices (all vertices before $u$ in $\pi$ are out-neighbors of $u$ and all vertices after $u$ are in-neighbors).
So, $u$ must receive a unique color, hence without loss of generality we have that $L_1 = \{u\}$. 

Next, consider vertex $v$ in the coloring of $L_\pi(G)$. Without loss of generality we have that $v \in L_2$. 
We claim that $|L_2| \le s$. Indeed, notice first that  $L_2$ is a clique in $R_\pi(G)$.
Who are the neighbors of $v$ in $R_\pi(G)$? Since $X$ appears before $v$ in $\pi$, no vertex of $X$ is a neighbor of $v$ in $R_\pi(G)$.
Since $Y$ appears after $v$ in $\pi$, no vertex of $Y$ is a neighbor of $v$ in $R_\pi(G)$.
Since $u$ appears after $v$ in $\pi$, and since $(u,v) \in E(G)$, $u$ is not a neighbor of $v$ in $R_\pi(G)$.
Consider some vertex $w \in N^-(u)$. Observe that $u$ is a neighbor of $v$ in $R_\pi(G)$ if and only if the $3$-set
$\{u,v,w\}$ forms a cyclic triangle. So, the set of neighbors of $v$ in $R_\pi(G)$ is precisely $Q$,
and they are all out-neighbors of $v$. But, in $R_\pi(G)$ the order of the vertices of $Q$ is according to $\sigma$,
so they induce no transitive sub-tournament of order $s$. So, together with $v$, there is no transitive sub-tournament of
order $s+1$ in the subgraph of $R_\pi(G)$ induced by $Q \cup \{v\}$. Thus, $|L_2| \le s$.

So, we have that
$$
n^2-\lceil (2-\frac{1}{\sqrt{2}})n \rceil + 3 = m = \sum_{i=1}^{n} |L_i| \le 1 + s + (n-2)n
$$
$$
= 2 + \lfloor \sqrt{2q} \rfloor + (n-2)n
$$
$$
\le 2 + \lfloor \sqrt{\frac{m+1}{2}} \rfloor + (n-2)n
$$
$$
\le 2 + \sqrt{\frac{m+1}{2}} + (n-2)n
$$
$$
\le 2 + \frac{n}{\sqrt{2}} + (n-2)n
$$
$$
= n^2- (2-\frac{1}{\sqrt{2}})n + 2
$$
which is a contradiction. Notice that in the fourth inequality we have used
$m  \le n^2-1$ which is true for $n \ge 3$ since $m=n^2-\lceil (2-\frac{1}{\sqrt{2}})n \rceil + 3$.
\qed

\section{Lower bound for $g(n)$}

\begin{lemma}\label{l:points}
Suppose $p_1=(x_1,y_1), \ldots, p_{2n}=(x_{2n},y_{2n})$ are $2n$ distinct points of $[n]^2$.
Then there are four points $p_i,p_j,p_k,p_l$ such that
$y_i=y_j < y_k = y_\ell$, $x_i < x_j, x_k < x_\ell$ and $[x_i,x_j] \cap [x_k,x_\ell]$ contains at least two integers.
\end{lemma}
\Proof
For $r=1,\ldots,n$ let $I_r$ be the consecutive set of integers whose smallest element is $\min\{x_i~|~ p_i=(x_i,r)\}$ and whose largest element is
$\max\{x_i~|~ p_i=(x_i,r)\}$. Notice that it is possible that $I_r = \emptyset$.
Now, since there are $2n$ distinct points, the sum of the sizes of the $I_r$ is at least $2n$, so there are two of them that intersect in at least two
integer points. \qed

\vspace*{10pt}
\noindent
Observe that lemma \ref{l:points} is tight. Indeed, the set of points $(i,i),(i+1,i)$ for $i=1,\ldots,n-1$ in addition to the point $(n,n)$ are $2n-1$
points in $[n]^2$ for which the lemma does not hold.

Let $G_n$ be the tournament on vertex set $[n]^2$ defined as follows.
$((i,j),(k,\ell))$ is an edge  if $i < k$ and $j \neq \ell$ or if $i > k$ and $j = \ell$ or if $i=k$ and $j < \ell$.
So, for every $i=1,\ldots,n$ there is a ``right to left'' transitive tournament on $(n,i),(n-1,i),\dots,(1,i)$ but for any permutation
$\sigma \in S_n$, the vertices $(1,\pi(1)),\dots,(n,\pi(n))$ form a ``left to right'' transitive tournament.

\vspace*{10pt}
\noindent
{\bf Proof of Theorem \ref{t:2}.}\,
We prove that $f(G_n) \le 3n^{7/4}$. This implies that
$g(n) \ge (n/3)^{8/7} \ge n^{8/7}/4$, as required.
To prove this fact, we will show that for every linear ordering $\pi$ of the vertices of $G_n$,
$\chi(L_\pi(G_n)) \le 3n^{7/4}$. This suffices since every spanning acyclic subgraph can be topologically sorted so that
it is a subgraph of $L_\pi(G_n)$ for some $\pi$.
Proving that $\chi(L_\pi(G_n)) \le 3n^{7/4}$ is equivalent to proving that $R_\pi(G_n)$ can be covered
by at most $3n^{7/4}$ transitive tournaments, since these tournaments are independent sets in $L_\pi(G_n)$.

We say that the vertices $(1,i),\dots,(n,i)$ are a {\em row} of $G_n$ and recall that each row induces a transitive tournament.
A {\em row clique} is a subset of vertices of a row. Thus, a row clique induces a transitive tournament.

So, consider some permutation $\pi$ and the corresponding $R_\pi(G_n)$.
We prove that $R_\pi(G_n)$ can be covered by at most $3n^{7/4}$ transitive tournaments.
We first find a maximum set of pairwise vertex-disjoint transitive tournaments in $R_\pi(G_n)$,
under the additional requirement that each such tournament has size at least $n^{1/4}/2$.
Assume that we have found $t$ such transitive tournaments, $X_1,\ldots,X_t$ (possibly $t=0$ if there are no tournaments meeting the size requirement).
Now, $|X_i| \ge n^{1/4}/2$ for $i=1,\ldots,t$ and $\sum_{i=1}^t |X_i| \ge tn^{1/4}/2$. In particular, $t \le 2n^{7/4}$.

Consider the remaining vertices not in these transitive tournaments, $Y = V(G_n) \setminus \cup_{i=1}^t X_i$.
If $|Y| \le n^{7/4}$ we are done.
Indeed, taking the elements of $Y$ as singletons, and taking the $t$ sets $X_1,\ldots,X_t$, we obtain a covering of
$R_\pi(G_n)$ with $t+|Y|$ transitive tournaments, and $t+|Y| \le 3n^{7/4}$, as required.

So we may now assume that the remaining set $Y$ has $|Y| \ge n^{7/4}$.
Let $Y_i$ be the elements of $Y$ that belong to the $i$'th row, and let $Y_i=\{(i_1,i),\ldots,(i_{r_i},i)\}$
where $i_1 < i_2 < \cdots < i_{r_i}$. Observe that $\sum_{i=1}^n r_i = |Y|$. More importantly, the
vertices of $Y_i$ as they appear in the linear ordering $\pi$, have no decreasing subsequence of order $n^{1/4}/2$ with respect to the
sequence $i_1,\ldots,i_{r_i}$ because such a decreasing sequence corresponds to a row clique in $R_\pi(G_n)$ of size at least $n^{1/4}/2$,
contradicting the maximality of $t$
(for example, if $\pi((7,i)) < \pi((13,i)) < \pi((9,i))$, then $(13,i),(9,i)$ form a decreasing sequence of length $2$ in $\pi$ and they induce a transitive tournament
of order $2$ in $R_\pi(G_n)$ since $((13,i),(9,i)) \in E(G_n)$).

From a theorem of Erd\H{o}s and Szekeres \cite{ES-1935} it now follows that each $Y_i$ has an increasing sequence in $\pi$, of order at least $2r_in^{-1/4}$.
So, the overall length of all these increasing sequences is at least $2|Y|n^{-1/4} \ge 2n^{3/2}$.
Let $S \subseteq Y$ be the union of these increasing sequences where $S_i \subseteq Y_i$ is the longest increasing sequence of
vertices in row $i$.

Define a directed vertex-labeled graph $R$ whose vertex set is $S$ and the label of some $(x,i) \in S$ is $\pi((x,i))$ (its position in $\pi$). 
There is a directed edge in $R$ from some $(x,i)$ to some $(y,j)$ if $x < y$, $i \neq j$, and $\pi((x,i)) < \pi((y,j))$.

We first show that $\chi(R) \ge n^{1/2}$. To this end, let us prove that $R$ has no independent set of size $2n$.
Since the elements of $R$ are vertices of $G_n$, we can think of any set of $2n$ vertices of $R$ as $2n$ points in $[n]^2$.
By Lemma \ref{l:points}, such a set has four points $(x_i,i),(x_j,i),(x_k,k), (x_\ell,k)$ such that
$i < k$, $x_i < x_j$, $x_k < x_\ell$ and $[x_i,x_j] \cap [x_k,x_\ell]$ contains at least two integer points.
Thus, the cases to consider are either (i) $[x_k,x_\ell] \subseteq [x_i,x_j]$ (ii) $[x_i,x_j] \subseteq [x_k,x_\ell]$
(iii) $x_i < x_k < x_j < x_\ell$ (iv) $x_k < x_i < x_\ell < x_j$.

Consider first case (i).
If $\pi((x_i,i)) < \pi((x_k,k))$ and $x_i < x_k$, then $(x_i,i),(x_k,k)$ are adjacent in $R$.
Otherwise, if $\pi((x_i,i)) > \pi((x_k,k))$ and $x_i < x_k$, then since $\pi((x_i,i)) < \pi((x_j,i))$ (recall that in $\pi$, $S_i$ is an increasing sequence)
we have that $\pi((x_k,k)) < \pi((x_j,i))$ and since $x_k < x_\ell \le x_j$ we have that $(x_k,k),(x_j,i)$ are adjacent in $R$.
Finally, if $x_i=x_k$ then if $\pi((x_i,i)) < \pi((x_k,k))$, then since $\pi((x_k,k)) < \pi((x_\ell,k))$ and since $x_i < x_\ell$
we have that $(x_i,i), (x_\ell,k)$ are adjacent in $R$. If, however, $\pi((x_i,i)) > \pi((x_k,k))$ then since $\pi((x_i,i)) < \pi((x_j,i))$
and since $x_k < x_j$ we have that $(x_k,k), (x_j,i)$ are adjacent in $R$.
The case (ii) is proved analogously to contain two adjacent vertices.

For case (iii) we have the following argument.
If $\pi((x_i,i)) < \pi((x_k,k))$, then $(x_i,i),(x_k,k)$ are adjacent in $R$.
Otherwise, $\pi((x_k,k)) < \pi((x_i,i)) < \pi((x_j,i))$ and since $x_k < x_j$ we have that $(x_k,k),(x_j,i)$ are adjacent in $R$.
The case (iv) is proved analogously to contain two adjacent vertices.

We have shown that $R$ has no independent set of size $2n$. Since $|V(R)|=|S|\ge 2n^{3/2}$ we have that $\chi(R) \ge 2n^{3/2}/(2n)=n^{1/2}$.
We would like to claim that $R$ has large transitive sets that intersect each row at most once. However, we cannot do that directly since
long paths in $R$ may alternate between only a few distinct rows. To overcome this, we partition the edge set of $R$ into two parts,
the upward edges $E_{up}$ and the downward edges $E_{down}$ as follows. Consider some edge $((x,i),(y,j))$ of $R$.
So we have $x < y$, $i \neq j$ and $\pi((x,i)) < \pi((y,j))$. We add this edge to $E_{up}$ if $i <j$ and to $E_{down}$ if $i > j$.
Let $R_{up}$ be the spanning subgraph of $R$ induced by the edges $E_{up}$ and let $R_{down}$ be the spanning subgraph of $R$ induced by the edges $E_{down}$.
By the product coloring argument mentioned in the introduction, we have that $\chi(R_{up}) \cdot \chi(R_{down}) \ge \chi(R)$.
Since $\chi(R) \ge n^{1/2}$ we have that at least one of them, say, $R_{up}$, has $\chi(R_{up}) \ge n^{1/4}$ (the proof in the case where
$R_{down} \ge n^{1/4}$ is analogous).

So, $R_{up}$ is an oriented graph with $\chi(R_{up}) \ge n^{1/4}$. By the well-known Gallai-Hasse-Roy-Vitaver Theorem \cite{gallai-1968,hasse-1965,roy-1967,vitaver-1962}, this means that $R_{up}$ has a directed path of order at least $n^{1/4}$. But this implies that $R_\pi(G_n)$ still has a transitive tournament
with $n^{1/4}$ all of whose vertices are in $Y$ (in fact, in $S$), contradicting the maximality of $t$.

\qed

\section{Random tournaments}

{\bf Proof of Proposition \ref{prop:1}.}\,
Let $G$ be a random tournament on vertex set $[n]$. We must prove that there exist positive constants $c$ and $C$ such that
$cn/\log n \le f(G) \le Cn/\log n$ with probability $1-o_n(1)$.

Recall the definition of $R_\pi(G)$ for $\pi \in S_n$ given Section 2.
Observe that for each fixed $\pi$, the fact that $G$ is a random tournament means that $R_\pi(G)$ is the (usual undirected)
random graph ${\cal G}(n,1/2)$. Indeed, each edge of $G$ is independently an edge of $R_\pi(G)$ with probability $1/2$.

Consider first the lower bound. For this case, we need only to consider one specific permutation $\pi$ (say, the identity).
It is well-known \cite{bollobas-1988} that the chromatic number of a random graph $H \in {\cal G}(n,1/2)$ is $\Theta(n/ \log n)$ with probability $o_n(1)$,
and in particular, there is a constant $c$ such that with probability $1-o_n(1)$, $\chi(H) \ge cn/\log n$
(in fact, any constant $c < 1/2$ suffices).
In other words, $\chi(R_\pi(G)) \ge cn/\log n$ with probability $1-o_n(1)$ and consequently
$f(G) \ge cn/\log n$ with probability $1-o_n(1)$.

As for the upper bound, slightly more care is needed. We would like to prove that there exists a positive constant $C$ such that
for any given permutation $\pi \in S_n$,
$$
\Pr[\chi(R_\pi(G)) \ge \frac{Cn}{\log n}] \le \frac{1}{(n+1)!}\;.
$$
Indeed, if we do that, then by the union bound applied to all $n!$ permutations, we will obtain
$$
\Pr[f(G) \ge \frac{Cn}{\log n}] \le \frac{1}{n+1} = o_n(1)\;.
$$
So, equivalently, our task is to prove that for a random graph $H \in {\cal G}(n,1/2)$, 
$$
\Pr[\chi(H) \ge \frac{Cn}{\log n}] \le \frac{1}{(n+1)!}\;.
$$
It follows from the proof of Bollob\'as \cite{bollobas-1988}, as well as the argument in \cite{AS-2004},
that for $G \in {\cal G}(n,1/2)$, if we set $m=\lfloor n/\ln^2 n \rfloor$, then the probability that some set of $m$ vertices
has a maximum independent set of size smaller than $2(1-o_n(1))\log n$ is at most $e^{-m^{2-o(1)}}$. So, in particular, for all $n$ sufficiently large, this latter probability is
less than $1/(n+1)!$. In other words, the following greedy procedure succeeds with probability at least $1-1/(n+1)!$.
Take an independent set of size at least $2(1-o_n(1))\log n$ and color it with color $1$, removing it from the graph. Continue to find another
independent set of this cardinality, coloring it with color $2$, removing it from the graph, and so on. Halt when the number of remaining vertices
is less than $m$ and color the remaining vertices each with a separate color. The resulting coloring uses only $(1+o_n(1))n/(2\log n)$ colors.
Hence, by the union bound,
$$
\Pr[\chi(H) \le (1+o_n(1))\frac{n}{2\log n}] \ge 1- e^{-m^{2-o(1)}} n \ge 1-\frac{1}{(n+1)!}\;.
$$
So, for any $C > 1/2$ we have that for all $n$ sufficiently large,
$$
\Pr[\chi(H) \ge \frac{Cn}{\log n}] \le \frac{1}{(n+1)!}\;.
$$
\qed

\section{$f(2)=3$}
In this section we prove Theorem \ref{t:3}.
Recall that a {\em subdivision} of an undirected graph $H$ is obtained by replacing each edge with
a path of length at least $1$. A subdivision is {\em odd} if each edge is replaced with a path
of odd length. We call a cycle in an oriented graph {\em inconsistent} if it is a cycle in the underlying graph but not a directed cycle in the oriented graph. We need the following two lemmas:
\begin{lemma}\label{l:51}
Every orientation of a non-bipartite $K_4$-subdivision has an inconsistent odd cycle.
\end{lemma}
The second lemma is a result of Thomassen \cite{thomassen-2001}, proving a conjecture of Toft \cite{toft-1975}.
\begin{lemma}\label{l:52}
Suppose $G$ is a $4$-chromatic graph, then $G$ contains an odd $K_4$-subdivision. \qed
\end{lemma}

\noindent
{\bf Proof of Theorem \ref{t:3}.}
Let $G$ be  an orientation of a $4$-chromatic graph. By Lemma \ref{l:51} and Lemma \ref{l:52},
$G$ has an inconsistent odd cycle. This cycle is an acyclic subgraph of $G$ with chromatic number
$3$. This implies that $f(2) \le 3$. To see that $f(2) \ge 3$ just consider any directed odd cycle.
Any acyclic subgraph of such a cycle has chromatic number at most $2$. \qed

\noindent
{\bf Proof of Lemma \ref{l:51}.}
Consider a non-bipartite $K_4$-subdivision. It has four vertices with degree $3$, denoted by $u,v,w,x$.
As there is an odd cycle, it either passes through all $u,v,w,x$ or through three of them.
But observe that in the former case, there must also be an odd cycle that passes through three of them.
So, without loss of generality, the cycle $C$ passing through $u,v,w$ (but not through $x$) is of odd length.
Now consider an orientation of this $K_4$-subdivision.
If $C$ is inconsistent, we are done. Otherwise $C$ is directed and without loss of generality its direction is from $u$ to $v$, from $v$ to $w$
and then from $w$ to $u$.

We need some notation for the lengths of the six paths connecting the pairs of vertices with degree $3$.
The path between $u$ and $v$ is of length $a$.
The path between $v$ and $w$ is of length $b$.
The path between $w$ and $u$ is of length $c$.
Thus, $a+b+c$ is odd as it is the length of $C$.
The path connecting $x$ to $u$ is of length $d$.
The path connecting $x$ to $v$ is of length $e$.
The path connecting $x$ to $w$ is of length $f$.

Since $a+b+c$ is odd, at least one of $d+e+a$, $e+f+b$, $f+d+c$ is odd.
Assume without loss of generality it is $d+e+a$. So if the cycle through $u,v,x$ is inconsistent, we are done.
Otherwise, it is a directed cycle, so the orientation of the path is from $x$ to $u$ and from $v$ to $x$.
Now, if $e+f+b$ is odd, then we are done since the cycle through $v,w,x$ is inconsistent. Hence we can
assume that $e+f+b$ is even and similarly we may assume that $f+d+c$ is even.

Now, if $f+d+a+b$ is even, then $(f+d+a+b)+(f+d+c)$ is even which implies that $a+b+c$ is even, a contradiction.
Thus, $f+d+a+b$ is odd. So the cycle passing through $w,x,u,v$ (in this order), whose length is $f+d+a+b$ is of odd length.
If it is inconsistent, we are done. Otherwise, the path between $w$ and $x$ is oriented from $w$ to $x$.

Similarly, if $c+a+e+f$ is even, then $(c+a+e+f)+(e+f+b)$ is even which implies that $a+b+c$ is even, a contradiction.
Thus, $c+a+e+f$ is odd. So the cycle passing through $w,u,v,x$ (in this order), whose length is $c+a+e+f$ is of odd length.
Since it is inconsistent, we are done. \qed

\bibliographystyle{plain}

\bibliography{references}

\end{document}